\documentclass[12pt,a4paper]{article}
\usepackage{amsmath,theorem}
\usepackage{epsfig}
\usepackage{amssymb}
\usepackage[latin1]{inputenc}
\usepackage[T1]{fontenc}
\newtheorem{lemma}{{\bf Lemma}}
\newtheorem{theorem}[lemma]{{\bf Theorem}}
\newtheorem{definition}[lemma]{{\bf Definition}}
\newcommand{\sep}{$\!${\rm.}\ }
\title{Small holding circles}
\author{Augustin F{\sc ruchard}}
\date{\empty}
\begin{document}
\maketitle
\begin{quote}
{\small\sl
Laboratoire MIA, EA3993, 
Universit\'e de Haute Alsace, \\
4, rue des Fr\`eres Lumi\`ere,
F-68093 Mulhouse cedex,
France\medskip\\
E-mail: {\tt Augustin.Fruchard@uha.fr}
\\
Tel. (+33) 389 33 66 37
\\
Fax  (+33) 389 33 66 53
}
\bigskip\\ 
{\bf Abstract:}
A circle $C$ {\sl holds }a convex body $K$ if  $C$ does not meet  the interior of $K$ and if there does not exist any euclidean displacement which moves $C$ as far as desired from $K$, avoiding  the interior of $K$.
The purpose of this note is to explore how small can be a holding circle.
In particular it is shown that the diameter of such a holding circle can be less than the width $w$ of the body but is always greater than $2w/3$.
\bigskip\\
{\bf Keywords:} holding circle, immobilization, width, diameter, convex body.
\bigskip\\
{\bf Mathematical Subject Classification (2010):} 52A15, 52A40, 52B10.
\end{quote}
\section{ Introduction}
The question of holding a convex body has been often considered in the literature.
In \cite{c}, Coxeter asked about the minimal total length of edges of a cage holding the unit ball.
Besicovitch \cite{b} and Valette \cite{v} investigated this question.
Analogous problems of caging are still widely studied, see e.g.\! \cite{pvs,rb,vv}.
Concerning circumscribing polyhedra, Besicovitch and Eggleton show in \cite{be} that the polyhedron with minimal total length of edges enclosing the unit ball is a cube.
In \cite{b1}, Besicovitch determines the minimal length of a net holding  the unit ball.

F. Caragiu asked whether convex bodies exist that can be held by a very simple instrument such as a circle.
T.~Zamfirescu gives the answer in  \cite{z}:
not only such convex bodies exist, but they form a huge majority.
More precisely, they form a subset with dense interior among all convex bodies, with respect to the Hausdorff-Pompeiu distance.
Whether this subset is open or not is unclear up to my knowledge.

In the whole article, $K$ denotes a convex body of $\mathbb{R}^3$ which admits holding circles;
to shorten we say {\sl a holdable convex body}.
Let $w$ denote the width of $K$ (i.e. the minimal distance between two parallel planes enclosing the body), $D$ the minimal diameter of a circumscribing cylinder and $d$ the minimal diameter of a holding circle.
We will also consider the supremum, denoted $\delta $, of diameters of holding circles
(in the article, the word ``diameter'' refers sometimes to a segment, sometimes to its length). 

Several quantities may measure how large or small are holding circles:
e.g.\! the ratios $\frac dD,\ \frac dw,\ \frac\delta D$ and $\frac\delta w$.
Of course we have $\frac dD\leq\frac dw\leq\frac\delta w$ and $\frac dD\leq\frac\delta D\leq\frac\delta w$.
In Section 3.7, we provide an example, due to T.~Zamfirescu, where all these ratios are as large as desired.
Section 3.8 contains examples with ratios $\frac dD$ and $\frac\delta D$ as small as desired.
Concerning the ratio $\frac dw$, our main result is the following.
\begin{theorem}\label{t2}\sep
\sl
{\rm(a)} For any holdable convex body $K$, we have $\frac dw>\frac23$.
\medskip\\
{\rm(b)}  The constant $\frac23$ is sharp: there are convex bodies with $\frac dw$ as close to $\frac23$ as desired.
\end{theorem}
The surprising Item (b) gives a negative answer to a question of Jo\"el Rouyer, who wondered if $d$ would always be greater than or equal to $w$.
See Figure 3 for a family of bodies satisfying (b).
The same result holds for the ratio $\frac\delta w$, see Section 3.6.
In Section 3 we also discuss the higher dimensional case (Section 3.1), the case of few-vertex polyhedra (3.2 and 3.3), the cases of ``asymptotic'' equality (3.5) and related topics in the literature (3.9).
\section{Proof of Theorem \ref{t2}}
\begin{figure}[ht]
\vspace{-.3cm}
\begin{center}
\epsfysize6.0cm\epsfbox{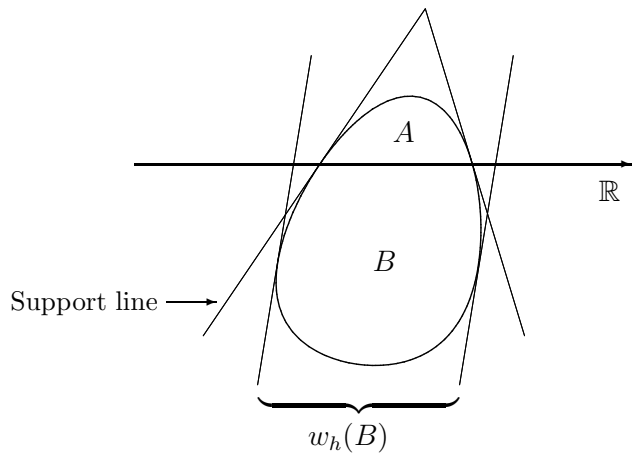}
\end{center}\vspace{-.6cm}
\caption{\label{fig1}
The {\sl subsets} $A,B$, two {\sl support lines} of $A\cup B$ intersecting $A\cap B$ and a {\sl strip} of minimal horizontal width containing $A\cup B$.}
\end{figure}

We first have to deal briefly with an aspect of planar convexity.
Given a non-horizontal strip in $\mathbb{R}^2$
$$
S=\{(x,y)\in\mathbb{R}^2\ ;\  ay+b_1\leq x\leq ay+b_2\}
$$
the {\sl horizontal width} of $S$ is $w_h(S)=b_2-b_1$.
Given a convex compact subset $B$ of $\mathbb{R}^2$, the {\sl horizontal width} of $B$,  $w_h(B)$, is the infimum of  $w_h(S)$ over all strips $S$ containing $B$, see Figure \ref{fig1}.
Of course, we have $w_h(B)\geq w_2(B)$, the usual planar width of $B$.
\begin{lemma}\label{l1}\sep
\sl 
Let $A,B$ be convex compact subsets of $\mathbb{R}^2$, $A$ in the closed upper half-plane  and $B$ in the closed lower half-plane, such that $A\cup B$ is also convex.
Then we have 
$$
w_h(A\cap B)=\min\big(w_h(A),w_h(B)\big)\mbox{ \rm and }w_h(A\cup B)=\max\big(w_h(A),w_h(B)\big).
$$
\end{lemma}

\noindent{\sl Proof of Lemma \ref{l1}} \sep
Observe that $A\cap B$ is a segment of $\mathbb{R}\times\{0\}$.
We obviously have $w_h(A\cap B)\leq\min\big(w_h(A),w_h(B)\big)$ and $w_h(A\cup B)\geq\max\big(w_h(A),$ $w_h(B)\big)$.
In order to prove that these are equalities, consider two support lines of $A\cup B$ at each end of $A\cap B$, see Figure \ref{fig1}.
If these lines can be chosen parallel, then we have $w_h(A\cap B)=w_h(A\cup B)$.
Otherwise they cross, say above, and then $A$ is included in a triangle of horizontal width $w_h(A\cap B)$.
This proves the first equality. 
For the second one, consider a strip of minimal horizontal length containing $B$.
The  boundary of this strip contains (at least) two points of the boundary of $B$ at the same altitude on each side, and which are not in $A\cap B$
(remind that we are in the case where there are no parallel support lines of $A\cup B$ at each end of $A\cap B$).
Therefore such a strip must contain $A$, yielding $w_h(B)=w_h(A\cup B)$.
\hfill\framebox[2mm]{\ }
\medskip

We now fix a holdable convex body $K$ and a holding circle $C$ of minimal diameter $d$ in a horizontal position.

Let $H$ denote the plane containing $C$.
This plane cuts $K$ in two convex bodies: $K^+$ above and $K^-$ below.
Let $K_0$ denote their intersection: $K_0=K^+\cap K^-=K\cap H$.
We call a {\sl horizontal slice} of $K^+$ the intersection of $K^+$ with a horizontal plane.

For any $\theta\in[0,\pi[$, let $V_\theta$ denote the vertical plane making an angle $\theta$ with $0x$.
Let $A_\theta$, resp. $B_\theta$  denote the orthogonal projection of $K^+$, resp. $K^-$ in  $V_\theta$.
Notice that the orthogonal projection of the holding circle, which is a segment 
of length $d$, must contain the segment $A_\theta\cap B_\theta$.
\begin{definition}{\label{d1}}\sep
\rm
With the above notation, we call $K$ an {\sl iceberg} if there exists a holding circle and an orientation of its axis such that $w_h(A_\theta)<w_h(B_\theta)$ for all $\theta\in[0,\pi[$.
\end{definition}
At a first glance it seems that icebergs cannot exist:
if $A_\theta$ is narrower than $B_\theta$ for any direction $\theta$, then it seems that the circle could be released through $A$. Nevertheless, as shows Figure \ref{fig5}, icebergs do exist.
In Section 3.3, we prove that tetrahedra and five-vertex polyhedra cannot be icebergs.

Notice that the width $w$ of $K$ satisfies
$w\leq w_2(A_\theta\cup B_\theta)\leq w_h(A_\theta\cup B_\theta)$
for all $\theta\in[0,\pi[$.
It follows from Lemma \ref{l1} that, if $K$ is an iceberg, then $w\leq w_h(B_\theta)$ for all $\theta\in[0,\pi[$.
Notice also that, if $w_h(B_\theta)<w_h(A_\theta)$ for all $\theta\in[0,\pi[$, then $K$ is an iceberg: it suffices to change the orientation of $C$.
\medskip

\noindent{\sl Proof of Theorem \ref{t2}} \sep

{\noindent\bf(a)} \ Two cases occur.
Firstly, if  $K$ is not an iceberg, this means that for some values of $\theta$ we have  $w_h(A_\theta)\leq w_h(B_\theta)$ and for some other ones we have $w_h(A_\theta)\geq w_h(B_\theta)$.
By continuity, there is a value $\theta_0\in[0,\pi[$ 
for which  $w_h(A_{\theta_0})=w_h(B_{\theta_0})$. Then we obtain from Lemma \ref{l1} 
\begin{equation}
\label{1}
w\leq w_h(A_{\theta_0}\cup B_{\theta_0})= w_h(A_{\theta_0}\cap B_{\theta_0})\leq d.
\end{equation}
By the way, this shows that, if $d<w$ then $K$ is necessarily an iceberg.
\medskip

\begin{figure}[ht]
\vspace{.3cm}
\begin{center}
\epsfxsize11cm\epsfbox{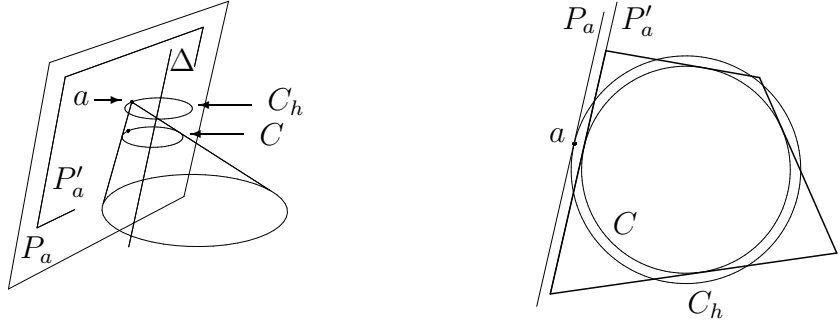}
\end{center}

\vspace{-.5cm}
\caption{\label{fig2} On the {\sl left}, the {\sl two circles} and the {\sl straight line} $\Delta$ joining their centers, the {\sl cone} with vertex $a\in C_h\cap K_h$ and the {\sl two planes} $P_a$ and $P'_a$. On the {\sl right}, the {\sl images} of $P_a$, $P'_a$, $C$ and $C_h$ by $\Pi$~;
in {\sl bold}, the {\sl boundary} of $I\cap H$.}
\end{figure}

\begin{figure}[ht]
\vspace{-1pt}
\begin{center}
\epsfxsize12cm\epsfbox{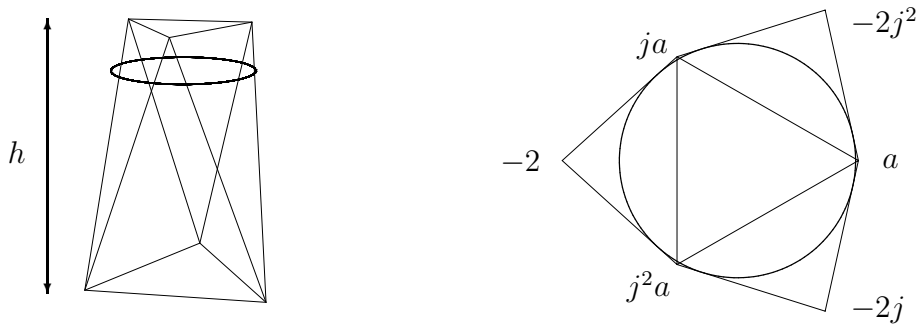}
\end{center}

\vspace{-15pt}
\caption{\label{fig5} an  octahedral {\sl iceberg} and its smallest {\sl holding circle}; here $a=1.38$ and $h=-5$.}
\end{figure}

Secondly, if  $K$  is an iceberg with $\frac dw<1$ (otherwise there is nothing to prove), then we have $d<w\leq w_h(B_\theta)$ for all $\theta$.
Since $C$ holds $K$, there is a horizontal slice $K_h$ of $K^+$ whose circumscribing circle $C_h$ has a diameter $d_h$ larger than $d$, otherwise $C$ could be released by a translation along the (continuous) curve of circumscribing centers of horizontal slices.
Let $\Delta$ denote the straight line joining the centers of $C$ and $C_h$. Let $\Pi$ denote the --- a priori non orthogonal --- projection of direction $\Delta$ into $H$
(we recall that $H$ is the horizontal plane containing $C$).
Let $\varphi$ denote the homothety of center the center of $C$ and of ratio $\frac d{d_h}$; in this manner, we have $C=\varphi(\Pi(C_h))$. 
Given $a\in C_h$, let $P_a$ denote the plane tangent to $C_h$ at $a$ and parallel to $\Delta$ and set $P'_a=\varphi(P_a)$; hence $P'_a$ is a plane parallel to $\Delta$ and tangent to $C$ at $\varphi(\Pi(a))$, see Figure \ref{fig2}.

For each $a\in C_h\cap K_h$, consider the cone of vertex $a$ and generatrix $C$. This cone contains $K^-$ because two points of $K^+$ and $K^-$ are joined by a segment which crosses the disk of boundary $C$. 
It follows that the closed half-space, denoted by $E_a$, containing $C$ and delimited by $P'_a$ contains $K^-\setminus K_0$ in its interior.
This means that $K^-\setminus K_0$ is in the interior of the intersection, denoted by $I$, of the half-spaces  $E_a$ for all $a\in C_h\cap K_h$.
Since $d<w_h(B_\theta)$ for all $\theta$, the width $w_h(B_\theta)$ is not attained close to the plane $H$ in the following sense:
there exists $\varepsilon =\varepsilon (\theta)>0$ such that $w_h(B_\theta)=w_h\big(B_\theta\cap\{z\leq-\varepsilon \}\big)$
(in fact, by compactness a single $\varepsilon $ is available for all $\theta$, but this is not necessary).
Since $K^-\cap\{z\leq-\varepsilon\}$ is a compact subset of  the interior of $I$, we obtain $w_h(B_\theta)<w_h(I_\theta)$ for all $\theta\in[0,\pi[$,
where $I_\theta$ denotes the orthogonal projection of $I$ into $V_\theta$.
The functions $\theta\mapsto w_h(B_\theta)$ and $\theta\mapsto w_h(I_\theta)$ are continuous on the compact set $[0,\pi]$, hence reach their infimum. 
Since $\Pi$ maps $I$ into $I\cap H$, we have
$\displaystyle\min_{\theta\in[0,\pi[}w_h(I_\theta)=w_2(I\cap H)$,
the planar width of $I\cap H$.

Because the circumscribing circle of $C_h\cap K_h$ is $C_h$ itself,
the circumscribing circle of the points $\varphi(\Pi(a)),\,a\in C_h\cap K_h$ is $C$ itself. Therefore $C$ is the greatest circle inscribed in $I\cap H$, hence satisfies $d\geq\frac23w_2(I\cap H)$, as is well-known since Blaschke \cite{bl}.
To sum up, we have
\begin{equation}\label{2}
w\leq\min_{\theta\in[0,\pi[}w_h(A_\theta\cup B_\theta)=\min_{\theta\in[0,\pi[}w_h(B_\theta)<\min_{\theta\in[0,\pi[}w_h(I_\theta)=w_2(I\cap H)\leq\frac{3d}2.
\end{equation}
{\noindent\bf(b)} \
With the identification $\mathbb{R}^3\simeq\mathbb{C}\times\mathbb{R}$ and the notation $j=\exp\big(\frac{2\pi i}3\big)$, choose $a>1$  (close to $1$) and $h>0$ (large), and consider $K$ the octahedron with vertices $(a,0)$, $(ja,0)$, $(j^2a,0)$, $(-2,-h)$, $(-2j,-h)$, $(-2j^2,-h)$. 
It has a horizontal holding circle $C$ of diameter $d =2a\cos\varphi$, where $\varphi=\arctan\big(\frac{a-1}{\sqrt3}\big)$, and of center $\big(0,-\frac{ha}{2\sqrt3}\sin(2\varphi)\big)$.
If $a$ tends to $1$ then $d$ tends to $2$ and if $h$ tends to infinity, then $w$ tends to $3$.
Hence the ratio $\frac dw$ is as close to $\frac23$ as desired.
Notice that the orthogonal projection in the horizontal plane shows $K$ as a hexagon close to an equilateral triangle and $C$ as its largest inscribed circle, see Figure \ref{fig5}. Observe also that three of the lateral faces are almost vertical.
\section{ Remarks and examples}
{\bf1.}
As suggested by the referee, Theorem \ref{t2} and its proof can be generalized in arbitrary dimension $n\geq3$:
if a convex body $K\subset\mathbb{R}^n$ of width $w$ admits a holding sphere of dimension $n-2$ and diameter $d$, then necessarily $\frac dw>C(n)$ where
$C(n)=\frac1{\sqrt{n-1}}$ if $n$ is even and $C(n)=\frac{\sqrt{n+1}}n$ if $n$ is odd.
The proof is the same, replacing the word ``circle'' by ``$(n-2)$-sphere'' and ``plane'' by ``hyperplane'', and using the Steinhagen inequality \cite{s}:
the width $w_{n-1}(K)$ and the inradius $r_{n-1}(K)$ of a convex body
$K\subset\mathbb{R}^{n-1}$ satisfy
$\frac{r_{n-1}(K)}{w_{n-1}(K)}\geq\frac{C(n)}2$,
see e.g. \cite{e} pp. 112--114 for a short proof.

To see that the constant $C(n)$ is sharp, 
we consider $\mathbb{R}^n$ euclidean with coordinates $x_1,\dots,x_n$ and split it in $\mathbb{R}^{n-1}\times\mathbb{R}$. 
In $\mathbb{R}^{n-1}$, consider the regular simplex, denoted $S_a$, centered at the origin and with a vertex at $(a,0,\dots,0)$, where $a>1$ is close to $1$.
Its side-length is $a\sqrt{\frac{2n}{n-1}}$.
Put this simplex in the hyperplane $x_n=0$; this is denoted $S_a\times\{0\}$.
With the same notation, consider the simplex $S_{1-n}$ and put it in the hyperplane $x_n=-h$ with $h>0$ large.
Then the convex hull of
$\big(S_a\times\{0\}\big)\cup\big(S_{1-n}\times\{-h\}\big)$
has a width close to $\frac2{C(n)}$ and a holding $(n-2)$-sphere of diameter $2a\Big(1+\frac{(a-1)^2}{n^2-2n}\Big)^{-1/2}$, hence close to $2$.
\medskip
\\
{\bf2.}
There is a short proof that $w\leq d$ for tetrahedra:
if $K$ is a polyhedron with a holding circle $C$ of minimal diameter $d$,
then $C$ contains at least two points of $K$ on two non-intersecting edges,
hence the diameter of $C$ is at least the distance between these edges,
i.e. the distance between the two parallel planes that contain each of these edges.
In the case of a tetrahedron, all the vertices, hence the whole tetrahedron, are in the closed spatial strip between these planes.
\medskip
\\
{\bf3.}
The example in Figure 2 has the minimal number of vertices required for an iceberg,
because neither tetrahedra nor five-vertex polyhedra can be icebergs.
Actually, assume by contradiction that a five-vertex polyhedron is an iceberg (the proof is similar for a tetrahedron).
With the notation above Definition \ref{d1}, among $K^+$ and $K^-$, one contains at most two vertices of $K$, say $K^+$.
Let the five vertices be labelled $a,b\in K^+$ and $c,e,f\in K^-$ and let $\theta_1$ be such that the vertical plane $V_{\theta_1}$ contains the direction of the edge $ab$.
Then in this direction we have $w_h(B_{\theta_1})\leq w_h(A_{\theta_1})$,
because otherwise the circle $C$ could be released by a translation along the axis joining its center to the middle of $a,b$. 
However, there is another direction $\theta_2$ such that $A_{\theta_2}$ is a triangle, yielding $w_h(B_{{\theta_2}})\geq w_h(A_{{\theta_2}})$: 
indeed if the vertices $a$ and $b$ are at the same altitude, then choose ${\theta_2}={\theta_1}+\frac\pi2\mod\pi$. 
Otherwise if $a$ is heigher than $b$, assume that, for the value $\theta=0$, the projection of $b$ on $V_\theta$ is, say, on the left of the polygon of the projections of $a,c,d,e$.
Then this projection of $b$ is on the right for the value $\theta=\pi$.
Therefore by continuity it has to cross this polygon for some $\theta_2\in[0,\pi]$. 
In conclusion, the inequalities $w_h(B_{\theta_1})\leq w_h(A_{\theta_1})$ and $w_h(B_{{\theta_2}})\geq w_h(A_{{\theta_2}})$ yield the contradiction.
\medskip
\\
{\bf4.}
For general holdable convex bodies that are not icebergs,
necessary conditions for equality $w=d$ can be derived from (\ref1):
the first equality  $w=w_h(A_{\theta_0}\cup B_{\theta_0})$ implies that the strip measuring $w_h(A_{\theta_0}\cup B_{\theta_0})$ has to be vertical;
secondly $K_0=K\cap H$ must contain all diameters of $C$ corresponding to the directions $\theta$ where $w_h(A_{\theta})=w_h(B_{\theta})=d$.

\begin{figure}[ht]
\vspace{2pt}
\begin{center}
\epsfxsize12.4cm\epsfbox{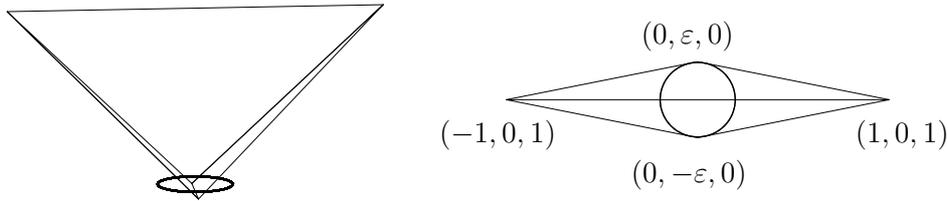}
\end{center}

\vspace{-12pt}
\caption{\label{fig4} a {\sl tetrahedron} with  small ratio $\frac dD$; here $\varepsilon =0.2$.}
\end{figure}

In particular, if $K$ is a tetrahedron such that $w=d$, then by projection in  $H$,
the four edges joining each vertex of $K^+$ to each vertex of $K^-$ form a rhombus:
they are tangent to $C$ and the points of tangency are on diameters, see e.g.\! Figure \ref{fig4}, right.
Because the distances between these points of tangency and the vertices of the rhombus are proportional to the distances between the plane $H$ containing $C$ and the vertices of the  tetrahedron,
it follows that the two other edges, one joining the vertices of $K^+$, and one joining those of $K^-$, are horizontal, orthogonal one to the other, and joined by their common orthogonal straight line in their middle.
To sum up, tetrahedra which satisfy $w=d$ are those with two non-intersecting orthogonal edges, joined by a common orthogonal straight line in their middle.
One can see that they form a $3$-dimensional submanifold of the $6$-dimensional space of congruence classes of tetrahedra in $\mathbb{R}^3$.
\medskip
\\
{\bf5.}
The case of ``asymptotic equality'' for Theorem \ref{t2} (a) can be described as follows.
For convenience, we use the framework of Nonstandard Analysis (NSA for short) but this is not essential:
instead of one nonstandard convex body $K$, the reader who is not acquainted with NSA may consider a whole sequence $(K_n)_{n\in\mathbb{N}}$.
Then expressions such as ``$a(K)$ is i-close to $b$'' (notation $a(K)\simeq b$), resp. ``$a(K)$ is i-large'' have to be replaced by ``there exists a subsequence  $(n_k)_{k\in\mathbb{N}}$ such that $a(K_{n_k})$ tends to $b$, resp. $a(K_{n_k})$ tends to $+\infty$, as $k\to+\infty$.

If $K$ is a convex body with $d=2$ and $w$ i-close to $3$, then all inequalities in (\ref2) have to be {\sl almost equalities}, i.e. equalities up to i-small numbers.
The last one $w_2(I\cap H)\simeq\frac{3d}2$ implies that $I\cap H$ is i-close to an equilateral triangle of height $3$ (with the Hausdorff-Pompeiu distance).
Here we use a well-known result due to Blaschke \cite{bl}:
given a planar convex set $A$ with planar width $w_2(A)$ and inradius $r(A)$, equality $w_2(A)=3r(A)$ holds only if $A$ is an equilateral triangle.
We now assume that $I\cap H$ is i-close to the triangle of vertices $(-2,0),(-2j,0)$ and $(-2j^2,0)$ where $j=\exp\big(\frac{2\pi i}3\big)$.
We can also describe the position of points of $K^+$ that prevent the holding circle to escape from the body.
Given any horizontal slice $K_h$ of $K^+$ with circumscribing circle $C_h$ of diameter $d_h>d$, we have that $d_h\simeq d$, that the segment joining the centers of $C$ and $C_h$ is almost vertical, and that all points of $K_h$ out of the horizontal circle of diameter $d$ and same center as $C_h$ must project on $C$ to points i-close to one of the three points  $(1,0),(j,0),(j^2,0)$. 
\medskip
\\
{\bf6.}
We now discuss other ratios than $\frac dw$.
Concerning the ratio $\frac\delta w$, one has the same result: this ratio, too, can be as close to $\frac23$ as desired.
To see this, we have to slightly modify the example of Figure \ref{fig5}, however, because this octahedron has also large holding circles. 
Actually, let $\alpha=\alpha(a),\beta=\beta(a)\in\,]0,1[$, $A=(z_A,0)$, $A'=(z_{A'},0)$ on two upper edges with 
$z_A=a\alpha+ja(1-\alpha),\,z_{A'}=a\alpha+j^2a(1-\alpha)$ 
and $B=(z_B,-h)$, $B'=(z_{B'},-h)$ on two lower edges with $z_B=-2\beta-2j(1-\beta),\,z_{B'}=-2\beta-2j^2(1-\beta)$
be such that $z_{A},z_{A'},z_{B},z_{B'}$
form a rectangle with diagonals orthogonal to the aforementioned edges
(i.e. such that ${\rm Im}\,z_A={\rm Im}\,z_{B'}$ and $\arg(z_A-z_B)=\frac\pi3$).
Then the segments $AB$ and $A'B'$ are diameters of a common holding circle. 
Nevertheless, it a possible to avoid these circles, either by moving slightly some vertices or by adding a seventh vertex, say $(0,1)$, in such a manner that there are no other holding circles than those close to $C$;
in particular holding circles of this seven-vertex polyhedron have diameters less than $2a$, yielding a ratio $\frac\delta w$ close to $\frac23$.
\medskip
\\
{\bf7.}
Concerning the ratio $\frac dD$, Tudor Zamfirescu presented in conferences the following ``bevelled cylinder'':
with $R>0$ arbitrarily large, let $x_1=(-R,-1,0)$, $x_2=(-R,1,0)$, $x_3=(R,0,-1)$ and $x_4=(R,0,1)$, let $C_1$ and $C_2$ be unit circles with axis the $0x$ axis, one centered at $(1-R,0,0)$ and the other at $(R-1,0,0)$;
then consider the convex hull of $\{x_1,x_2,x_3,x_4\}\cup C_1\cup C_2$.
One can see that any circle with diameter $[-R,R]\times\{0\}\times\{0\}$ which  does not cross the interior of the body (i.e.\! with axis far enough from the planes $y=0$ and $z=0$), is a holding circle of minimal diameter, hence $d=R$ whereas $D=1$.
\medskip
\\
{\bf8.}
A simple example with ratio $\frac dD$ as small as desired is the following one.
Given $\varepsilon >0$ arbitrarily small, the tetrahedron  with vertices 
$(0,\pm\varepsilon ,0)$ and $(\pm1,0,1)$ has a holding circle $C$ with diameter $d =2\sin\alpha$ where $\alpha$ is given by $\varepsilon =\tan\alpha$.
This is easily seen by orthogonal projection in a horizontal plane.
This circle is horizontal, centered on the $0z$ axis at altitude $d^2/4$, see Figure \ref{fig4}.

One can verify that the axis of the minimal circumscribing cylinder is the straight line $\big\{\big(x,0,\frac{1-\varepsilon ^2}2\big)\ ;\ x\in\mathbb{R}\big\}$ and its diameter is $D=1+\varepsilon ^2$, yielding $\frac dD$ arbitrarily small. 
Of course this tetrahedron has also large holding circles,
e.g.\! any circle of diameter $[(0,0,0),(0,0,1)]$ which does not cross the interior of the tetrahedron (i.e.\! with axis far enough from the planes $x=0$ and $y=0$),
but if we add the vertex $(0,0,-\varepsilon ^2)$ we obtain a five-vertex polyhedron with the same holding circle and with holding circles only close to $C$, thus with $\frac\delta D$ arbitrarily small. 

There exist also {\sl tetrahedra} with $\frac\delta D$ arbitrarily small,
e.g.\! the one with vertices
$A_1=(-2,-1,\varepsilon ),\ A_2=(-1,0,0),\ A_3=(2a,a,\varepsilon ),\ A_4=(a,0,0)$
with $\varepsilon >0$ arbitrarily small and $a=\varepsilon ^2$.
It is easy to verify that the circle in the plane $x=0$ of diameter $\varepsilon $ and center $\big(0,0,\frac\varepsilon 2\big)$ is actually  a holding circle.
Indeed, the edges $A_1A_3$ and $A_2A_4$ have the axis $0z$ as a common orthogonal straight line and the edges $A_1A_4$ and $A_2A_3$ cross the plane $x=0$ at
$\big(0,-\frac a{a+2},\frac{a\varepsilon }{a+2}\big)$
and $\big(0,\frac a{2a+1},\frac{\varepsilon }{2a+1}\big)$,
hence in the interior of the circle since $a<2\varepsilon^2$.
One can also verify that holding circles must be close to the origin and of small diameter.
\medskip
\\
{\bf 9.}
As a conclusion, we briefly describe related topics of the literature.
The first one is the problem of immobilization of convex bodies, a notion first introduced by W.~Kuperberg in \cite{k}.
In \cite{csu}, Czyzowicz, Stojmenovic and Urrutia prove that two-dimensional convex figures --- except circular disks --- can be immobilized by at most four points.
In \cite{bmu}, Bracho, Montejano and Urrutia show that three points suffice for convex figures bounded by a curve of class ${\cal C}^2$.
Mayer gives additional results and extensions in \cite{m}.
The three-dimensional case is studied by Bracho, Mayer, Fetter and Montejano in \cite{bfmm}:
a necessary condition for four points to immobilize a ${\cal C}^2$ convex body is that the four normal lines belong to one ruling of a quadratic surface.
These questions of immobilization are motivated by grasping problems in robotics, see e.g.  \cite{mnp,mnp1}.

A second related problem is to look for convex bodies passing through holes. Zindler \cite{zi} already considered an affine image of a cube passing through fairly small holes.
In \cite{iz}, Itoh and Zamfiresu look for the shape of a hole of minimal diameter and width through which can pass the regular unit tetrahedron $T$.
They find a hole of diameter $\sqrt3/2$, the width of a face of $T$, and of width $\sqrt2/2$, the width of $T$.
In  \cite{itz}, Itoh, Tanoue and Zamfiresu study the same tetrahedron passing through a circular and a square hole.
Triangular holes are considered by B\'ar\'any, Maehara and Tokushige in \cite{bmt} and higher dimensional holes in  \cite{mt}. 

Another related topic is known as Prince Rupert's problem, see  \cite{cfg}, problem B4. 
The original question is to cut a hole in the unit cube, large enough to let a larger cube passing through it. See also \cite{jw} for generalizations to rectangles.
\medskip\\
{\bf Acknowledgement:}
The author wishes to thank Tudor Zamfirescu warmly for his encouragements and for his highly valuable suggestions on a first version of the manuscript
and the referee for relevant comments on a second version.

\end{document}